\documentclass[11pt]{ncc}
\usepackage[utf8]{inputenc} 
\usepackage[russian]{babel}
\usepackage{nccmath,amsmath,amsfonts,amscd,amssymb,latexsym,longtable,bm,enumerate,textcomp}
\begin{document}
УДК 512.54

\begin{center} \textbf{ A. A. SHLEPKIN}
\end{center}

\begin{center} \textbf{PERIODIC GROUPS SATURATED WITH FINITE SIMPLE GROUPS OF LIE TYPE RANK 1}
\end{center}

It is said that the group $G$ {\it is saturated} with groups from the set of groups $\mathfrak{R}$ if any finite subgroup of $G$ is contained in a subgroup of $G$ that is isomorphic to some group in $\mathfrak{R}$ $[\ref{Shlepkin}].$
 In the Kourovka notebook $[\ref{kourovka}]$ the following question is posed 14.101:

\textit{Is it true that a periodic group saturated with finite simple groups of Lie type whose ranks
are bounded in the aggregate, is itself a simple group of Lie type?}

In this paper we give a partial answer to this question for periodic groups saturated by finite simple groups of Lie type  rank 1.

\textbf {THEOREM.} \textit {Let the periodic group $ G $ be saturated by finite simple groups of Lie type of rank $1.$
Then $G$ is isomorphic to a group of Lie type $1$
over a suitable locally finite field.}

\begin{center} \textbf{ А. А. ШЛЁПКИН}
\end{center}

\begin{center} \textbf{ПЕРИОДИЧЕСКИЕ ГРУППЫ, НАСЫЩЕННЫЕ КОНЕЧНЫМИ ПРОСТЫМИ ГРУППАМИ ЛИЕВА ТИПА РАНГА 1}
\end{center}

\begin{center}
 \textbf{ В в е д е н и е}
\end{center}

Говорят, что группа $G$ {\it насыщена} группами из множества групп $\mathfrak{R}$, если любая конечная подгруппа из $G$ содержится в подгруппе группы $G$, изоморфной некоторой группе из $\mathfrak{R}$ $[\ref{Shlepkin}]$. 
 В Коуровской тетради $[\ref{kourovka}]$ поставлен вопрос 14.101:

\textit{Верно ли, что периодическая группа, насыщенная конечными простыми группами  лиева типа, ранги которых 
ограничины в совокупности, сама являeтся простой группой лиева типа?}

В данной работе дается частичный ответ на этот вопрос для периодических групп, насыщенных конечными простыми группами лиева типа ранга 1.

\textbf{ТЕОРЕМА.} \textit{Пусть периодическая группа $G$ насыщена конечными простыми группами лиева типа ранга $1$.
Тогда $G$ изоморфна группе лиева типа ранга $1$ 
над подходящим локально конечным полем.}

\begin{center}
 \textbf{\S 1. Определения, известные факты}
\end{center}

\textbf{ОПРЕДЕЛЕНИЕ 1.} \textit{Пусть группа $G$ насыщена группами из множества групп $\mathfrak{R}$. Тогда мнжество  $\mathfrak{R}$ будем называть  насыщающим множеством для $G$}.

\medskip

\textbf{ОПРЕДЕЛЕНИЕ 2.} \textit{Пусть $G$ -- группа, $\mathfrak{X}$-- множество групп. Запись
$G\ \widetilde\in\ {\mathfrak{X}}$ означет, что $G$ изоморфна некоторой группе из $\mathfrak{X}$}.

\medskip

\textbf{ОПРЕДЕЛЕНИЕ 3.} \textit{Пусть $G$ --- группа, $K$ -- подгруппа $G,$ $\mathfrak{X}$ --- множество групп. Через
$$\mathfrak{X}_G(K)=\{H\mid K\leq H\leq G, H\ \widetilde\in\ \mathfrak{X}\}$$ 
будем обозначать множество всех подгрупп $H$ группы $G,$ содержащих подгруппу $K$ и изоморфных группам из множества $\mathfrak{X}.$ Если $1$ --- единичная подгруппа группы $G,$ то 
$$\mathfrak{X}_G(1)=\{H\mid H\leq G, H\ \widetilde\in\ \mathfrak{X}\}$$ 
будет обозначать множество всех подгрупп $H$ группы $G,$ изоморфных группам из множества $\mathfrak{X}$. Если из контекста ясно о какой группе $G$ идет речь, то вместо $\mathfrak{X}_G(K)$ будем писать  $\mathfrak{X}(K)$, и соответственно вместо  $\mathfrak{X}_G(1)$ будем писать $\mathfrak{X}(1)$.}

\medskip

\textbf{ПРЕДЛОЖЕНИЕ 1.} {\it Пусть $U=U_3(q)$, где $q$~--- cтепень простого нечетного числа.

$1$.\textit{ Если $q+1$  не делится на $4$, то силовская $2$-подгруппа из $U$  изоморфна полудиэдральной группе
$SD(m)=\langle a,b| a^{2^{m+1}}=b^2=1, a^b=a^{-1+2^{m}}\rangle$, где $2^m$ делит $q-1$, $2^{m+1}$ не делит $q-1$.}

$2$.\textit{ Если $q+1$ делится на $4$, то силовская $2$-подгруппа из $U$ изоморфна сплетенной группе
$W_r(m)=\langle a_1,a_2,b | a_{1}^{2^{m}}=a_2^{2^{m}}=b^2=e, a_1a_2=a_2a_1, a_{1}^{b}=a_2, a_{2}^{b}=a_1\rangle$,
где $2^m$ делит $q+1$, $2^{m+1}$ не делит $q+1$}}

[\ref{Alperin}, \ref{Alperin1}].

\medskip

\textbf{ПРЕДЛОЖЕНИЕ 2.} {\it Силовская $2$ - подгруппа из $L_2(r)$, где $r$~--- степень простого нечетного числа,
является группой диэдра}

 [\ref{Alperin1}].

\textbf{ПРЕДЛОЖЕНИЕ 3.}  {\it Пусть периодическая группа $G$
насыщена группами из множества конечных простых групп 
$$\{L_2(r), U_3(q)\mid r, q\ ~-\ \text{нечетные},\ r>3  \}.$$
Тогда 
 $$G\ \widetilde\in\ \{L_2(R), U_3(Q) | \ R, Q - \text{локально конечные поля нечетных характеристик}\}$$ }
[\ref{Lyt1}].

\textbf{ПРЕДЛОЖЕНИЕ 4.} {\it Пусть периодическая группа $G$ насыщена конечными простыми неабелевыми группами, и
в любой её конечной $2$-подгруппе $K$ все инволюции из $K$ лежат в центре $K$.
Тогда
 $$G\ \widetilde\in\ \{J_1, L_2(R),Re(F),U_3(Q),Sz(P)\mid R,F,Q,P- \text{локально конечные поля}\}$$}
[\ref{filsmj}].

\begin{center}
 \textbf{\S 2. Доказательство теоремы }
\end{center}

Предположим обратное, и пусть $G$ -- контрпример. Положим
 $$\mathfrak{M}=\{L_2(f), U_3(h), Sz(2^{2m+1}), Re(3^{2n+1})| f>3, h>2, m\geqslant 1,n\geqslant 1\}~-$$
 множество всех конечных простых групп лиева типа ранга 1,   
$$\mathfrak{D}=\{L_2(r), U_3(q)\mid r, q~-\textsl{нечетные},\ r>3, q\geqslant 3\},$$  
$$\mathfrak{C}=\{L_2(2^l), U_3(2^k), Sz(2^{2m+1}),Re(3^{2n+1})\mid l>2, k>1\}.$$
Тогда $\mathfrak{M}=\mathfrak{D}\cup\mathfrak{C}$ --- насыщающее множество для группы $G.$  Ясно, что $\mathfrak{M}(1)=\mathfrak{D}(1)\cup\mathfrak{C}(1)$.

\medskip

\textbf {ЛЕММА 1}. $\mathfrak{D}(1)\neq\varnothing$.

ДОКАЗАТЕЛЬСТВО. Предположим обратное. Тогда
  $\mathfrak M(1)=\mathfrak{C}(1)$, и взяв $\mathfrak{C}(1)$ в качестве насыщающего множества для $G$, получим по предложению 4, что 
  
$$G\ \widetilde\in\ \{L_2(R), U_3(Q), Sz(F), Re(P)\mid R, Q, F, P~-\text{локально конечные поля}\}.$$

Противоречие с тем, что $G$ --- контрпример.

 Лемма доказана.
 
 \medskip
 
 \textbf {ЛЕММА 2}. $\mathfrak{C}(1)\neq\varnothing$.
 
 ДОКАЗАТЕЛЬСТВО. Предположим обратное.
 Тогда $\mathfrak M(1)=\mathfrak{D}(1)$, и взяв $\mathfrak{D}(1)$ в качестве насыщающего множества для $G$,  по предложению 3 получим, что 
 $$G\ \widetilde\in\ \{L_2(R), U_3(Q)\mid R, Q~- \text{локально конечные  поля}\}.$$  Противоречие с тем, что $G$~--- контрпример.
 
 Лемма доказана.

\medskip

 В дальнейшем, через $S$ будем обозначать силовскую 2-подгруппу из $G.$
 
 \medskip
 
\textbf {ЛЕММА 3}. {\it $S$ --- бесконечная группа.}

ДОКАЗАТЕЛЬСТВО. Предположим обратное. По теореме Шункова [\ref{Lytkina3}, предложение 8] все силовские  2-подгруппы из $G$ конечны и  сопряжены c $S$. Следовательно,  $S$ содержит подгруппу, изоморфную силовской 2-подгруппе $S_X$ из $X\in\mathfrak{C}(1)$ (лемма 2), а также содержит подгруппу, изоморфную силовской 2-подгруппе $S_Y$ из $Y\in\mathfrak{D}(1)$ (лемма 1).
Пусть $S$ --- ранга 2. Тогда 
$$\mathfrak{C}(1)=\{X\mid X < G,\ X\simeq U_3(4)\}.$$
 Если  $S$ не изоморфна $S_X$, то
$S$ изоморфна  $S_Y$ для некоторго $Y\in\mathfrak{D}(1)$ и 
 либо $S$~--- сплетенная 2-группа, либо $S$~--- полудиэдральная группа (предложение 1). Но в обоих случаях $S$ не может содержать подгруппу, изоморфную $S_X$. Сдедовательно, $S\simeq S_X$.
В этом случае
все инволюции из $S$ лежат в $Z(S)$, и по предложению 2
$$\mathfrak{D}(1)=\{Y\mid Y < G,\ Y\simeq L_2(r)\ (r=3,5 (mod\ 8) \}.$$ 
Но тогда $G$ насыщена группами из множества
$$\mathfrak{M}(1)=\mathfrak{D}(1)\cup\mathfrak{C}(1)=\{Z\mid Z < G,\ Z\ \widetilde\in\ \{L_2(r)\ \{r=3,5 (mod\ 8), U_3(4)\}\}.$$
По предложению 4 $$G\ \widetilde\in\ \{L_2(R), U_3(4)\mid R~-\text{локально конечное поле}\}.$$  Изоморфизм $G\simeq L_2(R)$ невозможен по причине того, что $G$ содержит $X\simeq U_3(4)$, а $L_2(R)$ не содержит подгруппу, изоморфную $U_3(4).$ Следовательно, $G\simeq U_3(4).$ Противоречие с тем, что $G$ -- контрпример. 
 
  Итак, ранг $S$ больше 2. Пусть $A < S$, и $A$~--- элементарная абелева группа порядка 8. По условию насыщенности $S < K\in\mathfrak{C}(1)$. Поскольку в $L_2(r)$ и в $U_3(q)$ ($r, q$ -- нечетные) нет элементарных абелевых 2-подгрупп порядка 8 (предложения 1, 2), то по условию насыщенности $S < K\in\mathfrak{C}(1)$. Таким образом, все инволюции из $S$ лежат в $Z(S)$.  
 В этом случае $\mathfrak{M}(1)$ состоит из групп изоморфных группам из множества
 $$\{L_2(r)\ (r=3,5 (mod\ 8), L_2(2^l), U_3(2^k), Sz(2^{2m+1}), Re(3^{2n+1}\}.$$
По предложению 4  
 $$G\simeq\{L_2(R), U_3(Q), Sz(F), Re(P)\mid R,Q,F,P~- \text{локально конечные поля}\}.$$ 
 Противоречие с тем, что $G$ -- контрпример.

Лемма доказана.

\medskip

\textbf {ЛЕММА 4}. {\it Если $S$ содержит элементарную абелеву подгруппу порядка $8$, то  $S$ -- локально конечная  группа периода $4$, и все инволюции из $S$ лежат в $Z(S).$}

ДОКАЗАТЕЛЬСТВО. Пусть $A$~--- элементарная абелева подгруппа порядка 8 из группы $S$. По лемме 3 и $[\ref{KuzLyt1},\  $предложение 4$]$ $A\leq S_A\leq  S$, где $S_A$ -- бесконечная локально конечная подгруппа. Будем считать $S_A$ максимальной локально конечной подгрупой группы $S$, содержащей группу $A$. Из условия насыщенности вытекает, что $S_A$ периода не более 4, и все инволюции из $S_A$ лежат в $Z(S_A).$ Если все инволюции из $S$ лежат в $S_A$, то для любого  $s\in S$ группа $\langle A, s\rangle$ --- конечная 2-группа. По условию насыщенности $\langle A, s\rangle < R\in \mathfrak{C}(1)$. Следовательно, $|s|\le 4$, $S$ --- локально конечная группа (теорема Санова, [\ref{sanov}]), все инволюции из $S$ лежат в $Z(S)$ (условие насыщенности), и в этом случае лемма доказана. Пусть теперь $i$ --- инволюция из $S\setminus S_A. $ Возьмем инволюцию $j\in S_A.$ Ясно, что конечная группа $\langle i, j\rangle\cap S_A=D\neq 1$ не лежит в $S_A$.  В силу нормализаторного условия в конечных 2-группах можно выбрать такие элементы  
$ x\in N_S(D)\setminus S_A, y\in N_{S_A}(D)\setminus D, x^2\in D, y^2\in D$. Следовательно, конечная группа $K_1=\langle x, y, D\rangle\not < S_A$, $D < D_1=K_1\cap S_A$, $|D_1|\geqslant 2|D|.$
Продолжая этот процесс, мы всегда  сможем для любого наперед заданного положительного числа $m$ 
найти такую конечную подгруппу $K_m\not < S_A$, что $|K_m\cap S_A|>2^m|D|$. В силу строения $S_A$ (см. выше), начиная с некоторого $m$ все инволюции из  $D_m=K_m\cap S_A$ лежат в $Z(D_m)$
и образуют элементарную абелеву подгруппу порядка не менее 8. Возьмем теперь $y\in N_{K_m}(D_m)\setminus D_m$ и $y^2\in D_m$. Пусть $x$ -- произвольная инволюция из $S_A.$
Тогда $ \langle x, y, D_m\rangle$ -- конечная 2-группа, и из условия насыщенности получаем 
$xy=yx$, т. е. $y$ перестановочен со всеми инволюциями из $S_A$. Так как $S_A$ --- максимальная локально конечная подгруппа, содержащая $A$, то $y\in S_A$. Противоречие с выбором $y.$ Итак, все инволюции из $S$ лежат в $S_A$ и, как показано выше, в этом случае утверждение леммы имеет место.

Лемма доказана.

\medskip

\textbf {ЛЕММА 5}. {\it Если $S$ не содержит элементарных абелевых подгрупп порядка $8$, то $S$~--- черниковская группа ранга $2$.}

ДОКАЗАТЕЛЬСТВО. Действительно, в данном случае $S$ ограниченного ранга и по [\ref{shunkov}, \ref{lyt7}] $S$~--- черниковская. По условию леммы  ранг $S$ равен 2.

Лемма доказана.

Завершим  доказательство теоремы. Предположим, что группа $G$ содержит подгруппу $A$, и $A$ --- элементарная абелева подгруппа порядка 8. Пусть $\mathfrak{R}$--- множество силовских 2-подгрупп группы $G$, каждая из которых содержит подгруппу, изоморфную группе $A.$ Пусть $\mathfrak{T}$~--- множество всех остальных силовских 2-подгрупп группы $G$. Если одновременно  $\mathfrak{R}\neq\varnothing$ и  $\mathfrak{T}\neq\varnothing$, то  по [\ref{Lytkina2}, лемма 6] для любого наперед заданного натурального $m$,
найдутся такие  $R\in\mathfrak{R}$,  $T\in\mathfrak{T}$  со свойством
 $|T\cap R|>m$,
что невозможно по леммам 4,\ 5. 
 Следовательно, $\mathfrak{T}=\varnothing$, все инволюции из $S$ лежат в $Z(S)$ (лемма 4), и  
$\mathfrak M(1)$ состоит из групп, изоморфных группам из множества 
 $$\{L_2(r)\ (r=3,5 (mod\ 8), L_2(2^l), U_3(2^k), Sz(2^{2m+1}), Re(3^{2n+1}\}.$$
По предложению 4  
 $$G\ \widetilde\in\ \{L_2(R), U_3(F), Sz(Q), Re(P)\mid R,F,Q,P~-\text{локально конечные поля}\}.$$ 
Противоречие с тем, что $G$~--- контрпример. 
 Итак, группа $G$ не может содержать подгрупп изоморфных группе $A$.
Следовательно, $\mathfrak{R}=\varnothing$ и $\mathfrak{T}\neq\varnothing$. 
Таким образом, всегда $S\in\mathfrak{T}$. По леммам 3, 5 $S$~--- бесконечная черниковская  группа ранга 2. По лемме 2 $G$ содержит конечную подгруппу $K$ такую, что $K\in\mathfrak{C}(1).$ Пусть $S_K$~---
силовская 2-подгруппа из $K.$ Так как $\mathfrak{R}=\varnothing$, то $S_K$~--- ранга 2, а $K\simeq U_3(4).$ Поскольку $S_K$ лежит в некоторой силовской 2-подгруппе из $G$, а все силовские 2-подгруппы из $G$ бесконечные (лемма 3) и черниковские (лемма 5), то без ограничения общности можно считать, что $S_K\leq S.$ Возьмем $s\in S\setminus S_K$ такой, что $|s|>4$. По условию насыщенности $\langle s, S_K\rangle < F\in\mathfrak{D}(1)$.
  Ясно, что $F\simeq L_2(r),$ где $r$~--- нечетное. Обозначим через $S_F$ силовскую 2-подгруппу из $F$, содержащую $S_K.$ По предложению 2 $S_F~-$группа диэдра и не может содержать $S_K.$ Противоречие. 
   
Теорема доказана.

АДРЕС АВТОРА:

Шлепкин Алексей Анатольевич

Сиьирский Федеральный Университет

Г.Красноярск, пр. Свободный, 79

email: shlyopkin@gmail.com  
\end{document}